\documentclass[11pt]{article}
\usepackage[a4paper, top=1in, left=1in, right=1in, bottom=1in]{geometry}
\usepackage{setspace}
\usepackage{amsmath}
\usepackage{graphicx}
\usepackage{enumerate}
\usepackage{natbib}
\usepackage{xurl}
\usepackage{xcolor}
\usepackage{amssymb}
\usepackage[normal]{threeparttable}
\usepackage{booktabs}
\usepackage{hyperref}
\usepackage{subcaption}
\usepackage{tablefootnote}
\usepackage{url}

\begin{document}

\title{\bf Planning minimum regret $CO_2$ pipeline networks}
\author{Stephan Bogs$^a$, Ali Abdelshafy$^a$, Grit Walther$^a$ \\
$^a$Chair of Operations Management, RWTH Aachen University \\
Kackertstraße 7, 52072 Aachen, Germany}
\date{}
\maketitle
	
\begin{abstract}
The transition to a low-carbon economy necessitates effective carbon capture and storage (CCS) solutions, particularly for hard-to-abate sectors. Herein, pipeline networks are indispensable for cost-efficient $CO_2$ transportation over long distances. However, there is deep uncertainty regarding which industrial sectors will participate in such systems. This poses a significant challenge due to substantial investments as well as the lengthy planning and development timelines required for $CO_2$ pipeline projects, which are further constrained by limited upgrade options for already built infrastructure. The economies of scale inherent in pipeline construction exacerbate these challenges, leading to potential regret over earlier decisions.
While numerous models were developed to optimize the initial layout of pipeline infrastructure based on known demand, a gap exists in addressing the incremental development of infrastructure in conjunction with deep uncertainty. Hence, this paper introduces a novel optimization model for $CO_2$ pipeline infrastructure development, minimizing regret as its objective function and incorporating various upgrade options, such as looping and pressure increases. The model’s effectiveness is also demonstrated by presenting a comprehensive case study of Germany’s cement and lime industries. The developed approach quantitatively illustrates the trade-off between different options, which can help in deriving effective strategies for $CO_2$ infrastructure development.
\end{abstract}
		
\noindent%

\section{Introduction}
\label{sec:intro}

Carbon Capture and Storage (CCS) represents a pivotal technique in mitigating climate change by capturing $CO_2$ emissions directly from industrial sources and sequestering them in suitable geological formations. CCS will most likely be used by certain industries such as cement and lime, which are classified as hard-to-abate sectors due to their unavoidable process emissions \citep{abdelshafy_role_2022}. Therefore, CCS has been integrated as a key technique in several industrial transformation roadmaps \citep{iea_technology_2018}. The European Commission also introduced in 2023 the Net Zero Industry Act \citep{european_commission_net_2023}, which proposes an EU wide goal of providing $CO_2$ injection capacity at around 50 million tonnes by 2030. Germany recently outlined the key points for CCS development in Germany in their Carbon Management Strategy \citep{bmwk_eckpunkte_2024}.

To make CCS economically viable there arises the need to transport huge quantities of $CO_2$ over long distances from the emitters to potential sinks. Therefore, IEA outlines $CO_2$ transport and storage as one of the key working fields \citep{iea_co2_2023}. \citet{abdelshafy_role_2022} show that pipelines will play a crucial role as it is the most cost-effective transportation mode comparing it to trucks or railway. To reach the ambitious climate goals, the planning and installation of the pipeline infrastructure needs to start now. However, pipelines also come with high initial investment costs as well as ongoing operational and maintenance costs. These high costs go along with other challenges. Besides the associated techno-economic risks, there are several uncertainties regarding carbon prices, political support and social opposition \citep{abdelshafy_role_2022}. Also, there are still some controversies regarding the CCS demand of some sectors (e.g. steel) as there are other alternative decarbonization technologies (e.g. hydrogen direct reduction) \citep{material_economics_industrial_2019, iea_iron_2020}. Additionally, due to the immense size of the prospective networks, they are going to be established incrementally. For example, the studies of \citet{morbee_optimal_2011}  and \citet{morbee_optimised_2012} show that the CO2 network in EUrope will expand gradually over the next three decades. Therefore, based on the effectiveness of the planning approach, each construction phase can impact the next one either positively or negatively.

These ambiguities are specifically hindering fast developments using current planning methods. Decision makers need approaches that allow them to secure their investment against the uncertain environment. Herein, no-regret and low-regret are two terms often used in the discussion \citep{zep_trans-european_2020, agora_energiewende_no-regret_2021, in4climatenrw_ed_co_2021}. However, this concept is so far mainly addressed from a qualitative perspective. But instead of dealing with regret as an abstract concept, there is a genuine need to quantify and minimize it, thereby supporting the decision-making process. This poses the question of whether there are planning methods, that lead to minimum regret for $CO_2$ pipeline planning under uncertain demand. In this regard, demand refers to certain industries or countries joining the system at different times. 

Here, the main challenge is the economic features of pipeline construction projects. As pipeline construction has nonlinear cost functions and constraints \citep{parker_using_2003, benrath_co2_2020}, economies of scale play a major role in minimizing the construction costs. Also, once the pipeline is built, it cannot be upgraded to a bigger diameter later, without incurring extremely high costs. The application-oriented literature (e.g. \citet{mischner_gas2energynet_2015}) contains several techno-economic aspects of enhancing capacity within natural gas networks. This may involve methods such as adding more pumping stations and compressors along the existing pipeline to increase the pressure, or building parallel lines. However, the expenses associated with such measures remain significantly higher compared to initially building the pipeline with a high capacity. Thus, it might be advantages to build bigger pipelines in earlier phases even if they are not yet fully utilized. Nonetheless, oversizing poses risks as the potential emitters may decide not to join the network at a later stage. In such case, oversized pipelines would incur unnecessary costs in construction, operation and maintenance.

Given the high uncertainty surrounding which emitters will join the system, it is essential to optimize for multiple scenarios within a single model. This necessitates deriving a regret-minimizing, multi-period model that addresses non-linear cost functions and constraints while considering only specific upgrade options in subsequent periods. The planning should also account for geographic features and model numerous future demand scenarios. This paper proposes such a model and demonstrates its efficacy via a case study.
In terms of paper structure, the following section provides an overview of related research and state how our contribution fits into it. Afterwards, we describe our methodology for designing in Section 3 and layout the formal model in Section 4. In Section 5 we present a case study to demonstrate the model and discuss the results. Finally, we conclude the paper with Section 6, highlighting the main outcomes and an outlook for future research.

\section{Related research}
\label{sec:related}

Due to the crucial role of CCS, the relevant supply chains and infrastructure systems have emerged as important themes in literature. \citet{tapia_review_2018} provide a comprehensive review of decision support systems for CCS. The study demonstrates the wide range of methods that have been developed to address the different operations along the CCS supply chain, i.e. $CO_2$ capture, transportation and storage. For pipeline networks design, there are also various approaches. There is an abundance of literature for natural gas as well as for the underlying graph problems, e.g. the Minimum Fixed Charge Network Flow Problem (FCNFP). However, we focus this section specifically on strategic planning methods of mid-to-long distance $CO_2$ pipeline networks. Based on the modeling requirements mentioned in the previous section, we split the literature analysis into four parts: the general design and scalability, the approach to geographic accuracy, considerations of uncertainty and how our contribution fits into the literature.

\subsection{General design and scalability}

Pipeline design and operation involve non-linear constraints, which make too detailed models hard to manage and tract. Therefore, all subsequent publications incorporate assumptions and simplifications to address this challenge. In order to build a pipeline network, the planner has to select appropriate pipeline sizes. Thus, a common approach is to model discrete pipeline sizes with an assigned flow as a Mixed Integer Linear Program (MILP) \citep{middleton_scalable_2009, sun_development_2017, bennaes_modeling_2024}. The MILP-based models can be either single period or multi-period models. The multi-period models (e.g. \citet{jones_designing_2022}) incorporate a stepped approach of building the network incrementally. \citet{middleton_simccs_2020} improved their model with piece-wise linear approximation, aiming at making the optimization more scaleable. This advancement has been further expanded upon by \citet{jones_designing_2022}. 
\citet{whitman_scalable_2022} also demonstrate that modeling CCS infrastructure design as an FCNFP is feasible. Since this poses an NP-hard problem, they provide several heuristics for it. Some earlier models (e.g. \citet{elahi_multi-period_2014}) adopt a simplified approach with a combination of fixed and linear cost. Following \citet{parker_using_2003}, a quadratic cost function based on the inner pipeline diameter can be used as an approximation of the discrete sizes. Herein, a squared cost function can be employed to determine the $CO_2$ flow. \citet{yeates_heuristic_2021} model the problem as a minimum cost flow problem and provide a comparison of different available heuristics as well as their own heuristic. 

\subsection{Geographic accuracy}
\label{subsec:geographic}

Geographic features such as slopes, population density, and existing infrastructure must be considered. \citet{onyebuchi_systematic_2018} conclude that $CO_2$ is neither toxic nor explosive, but it could still pose a risk for humans. Moreover, building $CO_2$ infrastructure close to residential areas may face opposition from the local community. Elevation profiles also pose significant challenges not just for building the pipeline but also for technical limitations on how $CO_2$ is transported through pipelines. Therefore, integrating geodata during the design phase of the pipeline networks is essential. Herein, there are various approaches with different accuracy levels to account for the geographic factors. The conventional method is assuming straight lines between the emitters \citep{elahi_multi-period_2014, sun_development_2017}. \citet{bennaes_modeling_2024} have advanced this approach by using offshore and onshore detour factors. \citet{middleton_scalable_2009} consider geographic features like slopes, waterways, national parks and state parks, as well as infrastructure like highways and railroads. They also make a distinction between wetlands and urban areas. \citet{yeates_industrial_2024} also adopt the same approach and incorporate more factors. They integrate population density and consider existing natural gas pipelines as positive factor, based on the assumption that construction is cheaper to build where pipeline infrastructure already exists.

\subsection{Uncertainty}
\label{subsec:uncertainty}

Uncertainty has also emerged as a crucial aspect in the relevant models. While many studies have typically assumed known demands, various models consider variables such as uncertain sink capacity \citep{middleton_cost_2018}. For the specific use case of enhance oil recovery (EOR), several multistage, stochastic models are proposed (\citep{tapia_optimal_2016, abdoli_novel_2023}. The study of \citet{damore_european_2019} also incorporates uncertain storage capacities and emission targets. \citet{han_multiperiod_2012} consider uncertain amounts of $CO_2$ emissions, product prices, and operating costs in their multistage stochastic model. There are also multiple models that account for fluctuating carbon prices, either by using a multistage stochastic model \citep{elahi_multi-stage_2017} or qualitative analysis \citep{sun_impact_2022}. Demand uncertainty was not covered sufficiently, except for some publications like \citet{nicolle_build_2023}. The latter paper also introduce the concept of regret. However, while it shows that building for certain demand and then upgrading later leads to maximizing of regret over multiple periods, they do not look at the network layout, but rather general pipeline construction.

\begin{table}[h]
	\centering
	\begin{threeparttable}			
		\caption{Comparison of selected models}	
		\label{tab:bib}		
		\begin{small}			
			\begin{tabular}{rccccc} 					
				\toprule							
				& Model & Geodata \& & Multi-Period & Uncertainty & Upgrade  \& \\
				&  & Routing &  &  & Regret \\
				\midrule
				\citet{middleton_simccs_2020} & MILP & \checkmark & $\times$ & $\times$ & $\times$  \\
				\citet{sun_development_2017} & MILP & $\times$ & \checkmark	& $\times$ & $\times$\\
				\citet{whitman_scalable_2022} & Heuristic & \checkmark & $\times$
				& $\times$ & $\times$ \\	
				\citet{yeates_industrial_2024} & Heuristic & \checkmark & $\times$ & $\times$ &  $\times$ \\			
				\citet{bennaes_modeling_2024} & MILP & \checkmark* & 
				$\times$ & $\times$  &  $\times$ \\
				\citet{middleton_cost_2018} & MILP & \checkmark & $\times$ & Sinks & $\times$  \\		
				\citet{abdoli_novel_2023} & MILP & $\times$ & \checkmark
				& Sinks &  $\times$ \\				
				This model & MILP & \checkmark & \checkmark & Emitters &\checkmark\\
				\bottomrule	
			\end{tabular}				
		\end{small}
		\begin{tablenotes}				
			\item {\footnotesize * Using different detour factors for overland and sea pipeline}			
		\end{tablenotes}
	\end{threeparttable}
	\centering
\end{table}

\vspace{0.1cm}

\subsection{Knowledge gaps \& Our contribution}
\label{subsec:contribution}

While significant research has been conducted on the aforementioned aspects, pipeline network models still require more upgrades and integration. As shown in Table \ref{tab:bib}, most of the enhancements have focused on the geodata and routing, followed by considering multi-period optimization. Additionally, there is no one study that addresses these aspects combined. Also, stochastic models have been the primary approach for addressing uncertainties such as storage capacities and $CO_2$ prices. Herein, there is a research gap for worse-case optimization considering uncertain emitter scenarios. Furthermore, \textit{min-max} regret as a target function as well as the use of upgrade functions like in \citet{mischner_gas2energynet_2015} were, to our best knowledge, never used in the context of $CO_2$ pipeline network planning. To close this gap, we propose a novel two-stage multi-scenario optimization model that uses \textit{min-max} regret as a target function. The derived approach also aims for comprehensive analyses by considering all relevant aspects listed in Table \ref{tab:bib}. We are building on the general ideas of \citet{middleton_simccs_2020} and \citet{yeates_industrial_2024} for geographic data as well as using piece-wise linear approximation. Similar to \citet{jones_designing_2022} we follow an approach to develop the network over multiple periods. We will further make use of upgrade functions, inspired by \citet{mischner_gas2energynet_2015}.

\section{Method}
\label{sec:method}

We have a set of scenarios $\mathcal{S}$. Each scenario $s \in \mathcal{S}$ has a set of nodes $\mathcal{V'}$ that are $CO_2$ emitters (sources) or sinks with an amount of emissions (emitter) or maximum capacity (sinks). We select an initial scenario $s_1 \in \mathcal{S}$ to start building the pipelines ($t=0$). After a period of $n_1$ years we take another set of investment decisions  ($t=1$). After this stage we still consider running the pipeline network for the next $n_2$ years.

For our optimization model we consider geographical accuracy to derive a graph and a cost function for building and maintaining a pipeline section.

To connect the nodes $\mathcal{V'}$ in a geographical accurate way we follow the approach to incorporate geographic data outlined by \citet{middleton_scalable_2009}, that was improved and applied for Germany by \citet{yeates_industrial_2024}. 

In short, we create a rasta-map with squares with a side length of 1.5km. Then for geographical features like terrain, rivers, existing infrastructure, population density or national parks we apply multipliers to the squares, creating a penalty on each square. From one node to another we compute the least penalized way through the raster map. We then reduce the found routes to a common graph by introducing transport nodes where the routes merge. A detailed description of the approach can be found in \ref{app:graph}. 

\begin{figure}[h]
	\centering
	\includegraphics[width=0.65\linewidth]{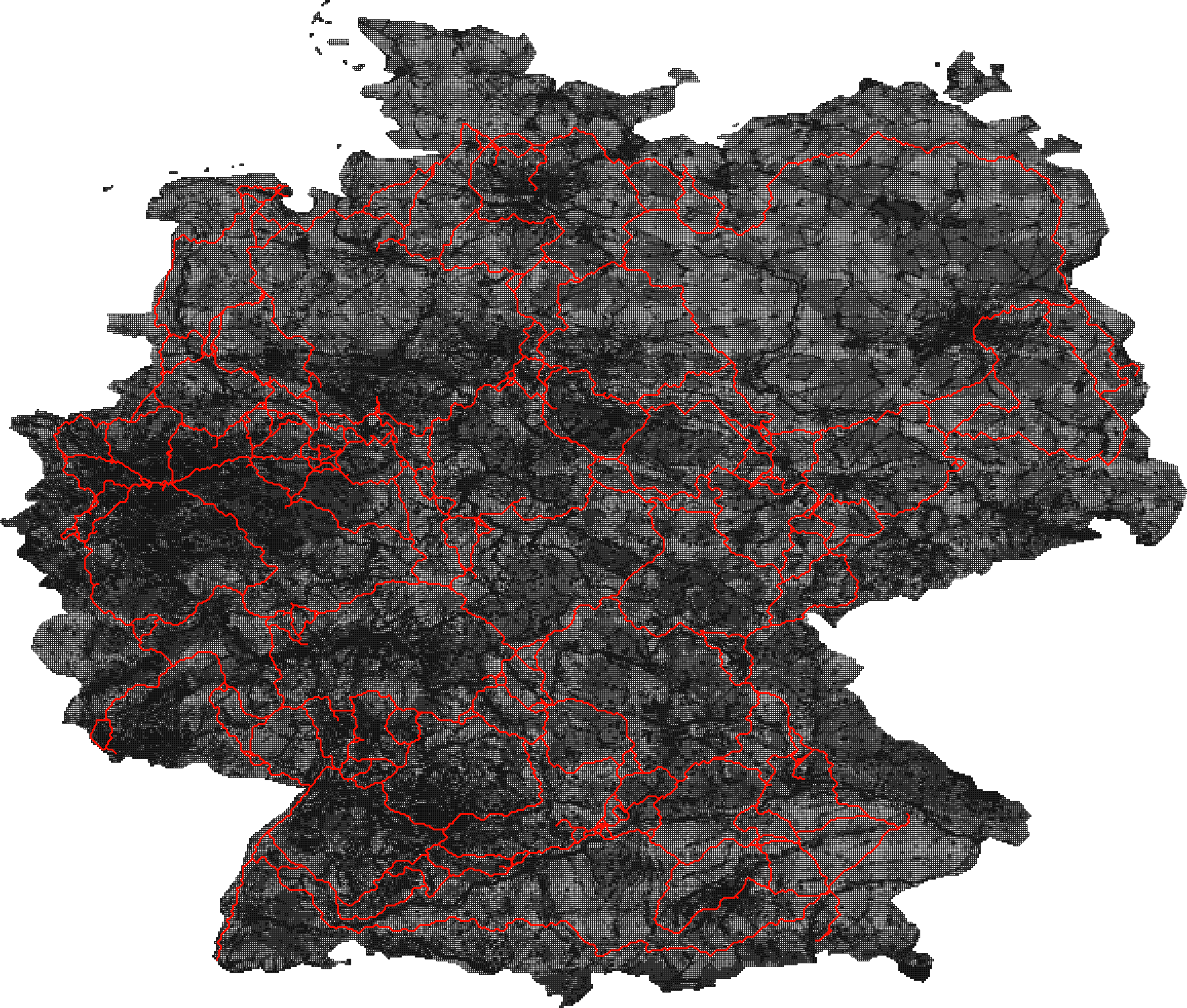}
	\caption{Rasta-map, lighter=better conditions, darker=worse conditions and  network graph (red)}
	\label{fig:rastamap}
\end{figure}

This results in a directed graph $\mathcal{G} = (\mathcal{V},\mathcal{A})$, where $\mathcal{V}'$ is a subset of $\mathcal{V}$, which is extended by the transport nodes. For an example see Figure \ref{fig:rastamap}. Furthermore for each arc $(i,j) \in \mathcal{A}$ we get the length $l_{ij}$ following the exact geographic line.

Each node $i \in \mathcal{V}$ has a demand $d_{is}$ depending on the scenario $s$. All transport nodes have $d_{is}=0$. Sinks have negative demands $d_{is} < 0$ and sources have positive demands  $d_{is} > 0$.
Sources that are not present in another scenario are transport nodes in that scenario. 

For costs we are using the cost function given by \citet{yeates_industrial_2024}, which is a slightly simplified version of \citet{parker_using_2003}: 

\begin{align*}
	cost(D) = (C1 \cdot D^2 + C2 \cdot D + C3) \cdot l_{ij}\\
\end{align*} 

$C1$, $C2$, $C3$ are cost constants and $D$ is the inner diameter of the pipe. To get this inner diameter, we use the following function of \cite{benrath_co2_2020}:

\begin{align*}
	D =\sqrt{\dfrac{F}{v \cdot \pi \cdot 0.25 \cdot p}}\\
	\text{$D$  inner diameter, $F$ flowrate, $v$ velocity, $p$ density}
\end{align*} 

Following the approach of \citet{yeates_industrial_2024} we fix pressure and temperature of the gas. We derive the density of the gas $p$, and we also fix velocity $v$ inside the pipeline. We now can derive a maximum flow for a given pipeline diameter $D$, and so derive the cost function based on a given flowrate $F$:

\begin{align*}
	cost(F) =  \left( \frac{F}{v \cdot \pi \cdot 0.25 \cdot p} \cdot C1 + \sqrt{\frac{F}{v \cdot \pi \cdot 0.25 \cdot p}} \cdot C2 + C3 \right) \cdot l_{ij}
\end{align*} 

\begin{figure}[h]
	\centering
	\includegraphics[width=0.7\linewidth]{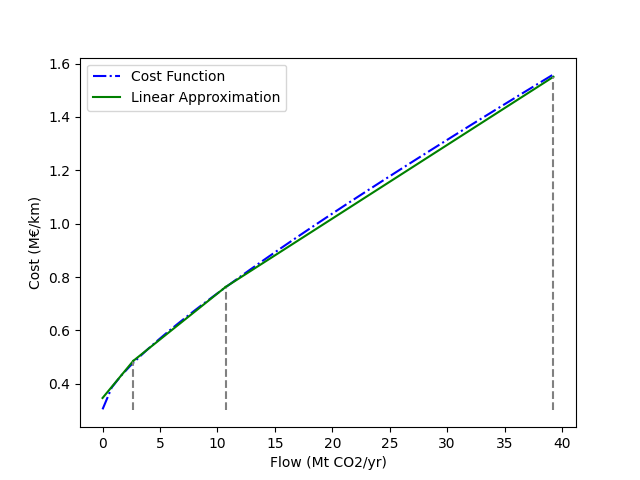}
	\caption[Cost function and approximation]{Cost function and approximation, 3 trends}
	\label{fig:costfunction}
\end{figure}

As this function is non-linear we approximate it with a set of piecewise linear functions, called trends $\mathcal{C}$. For each trend we derive a linear part $m_c$ and a fixed part $y_c$. An example how this approximation looks like can be seen in Figure \ref{fig:costfunction}. Yearly operations and management cost are taken as a percentage of the initial investment costs as in \cite{benrath_co2_2020}. Future cost are discounted using a discount rate $\tau$.

Based on the cost function and the graph we can compute a network in $t=0$: For each arc $(i,j) \in \mathcal{A}$  we have to decide if to build a pipeline with a specific cost trend $c \in \mathcal{C}$ and if so, how much maximum gas flow we have on that arc.

The following decisions have to be taken in $t=1$, e.g. five years after $t=0$: How do we extend the already given network in a cost optimal way. For this we have a set of upgrade operation $\mathcal{O}$, that we derived from \citet{mischner_gas2energynet_2015}. This include building parallel lines ($o=1$) (also called looping), and a pressure increase on an existing line ($o=2$). For the pressure increase we allow an extra flow, taking the original flow in $t=0$ multiplied by $o_2^{max}$. This comes at a cost factor of $o_2^{cost}$ times the original price.

\begin{figure}[h]
	\centering
	\includegraphics[width=0.9\linewidth]{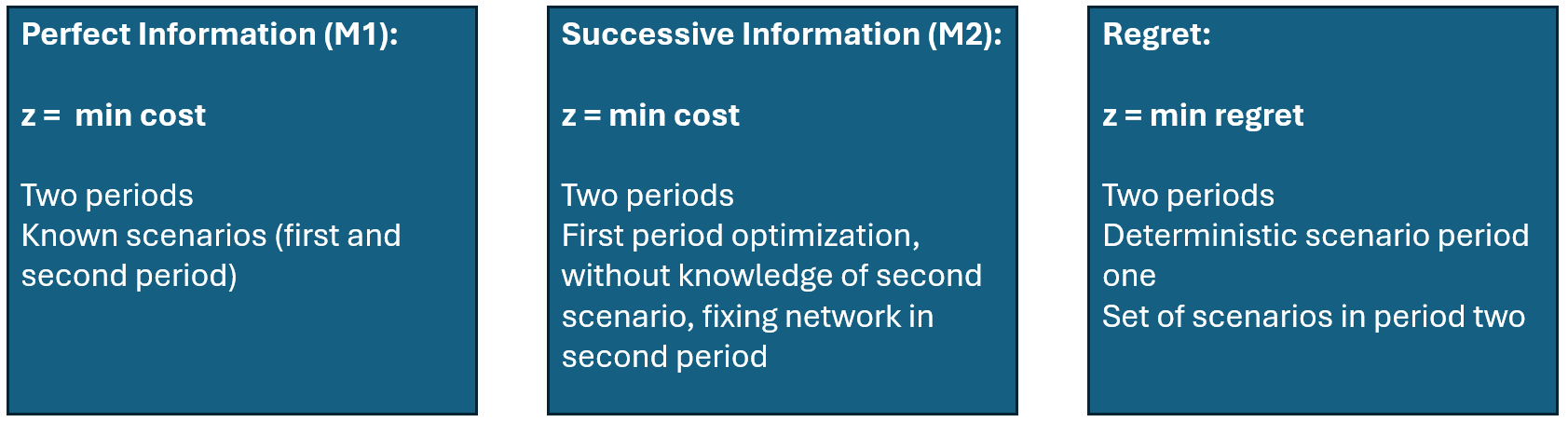}
	\caption{Overview over the different optimization models}
	\label{fig:models}
\end{figure}

Our main focus is a regret model, however to determine its efficacy we need two auxiliary models. Figure \ref{fig:models} gives an overview of them.

The first model (M1) is assuming perfect information of the scenario $s_1$ and $s_2$ in a phased approach. This leads to a layout in $t=0$ that can be easily extended in $t=1$. The successive information model (M2) computes an optimal network for a scenario $s_1$ without anticipation of $s_2$, which can lead to heavy interventions and regret of the decisions. 

For the regret model we then assume that all scenarios can become active in the next investment period $t=1$. We define a $min$ $max$ regret target function that considers the costs of all possible outcomes and computes the regret against the best solution, if the scenario outcome would have been known (M1). This gives the solution for $t=0$ that can be build, knowing the regret over all possible realizations in $t=1$ is minimal.

Our model is tractable compared to other publications as we reduce the economic model significantly. We do not include sink opening costs, sink sequestration cost, or investment cost at the emitter side. Each scenario assumes that all participants need to handle all $CO_2$. The focus is on the pipeline construction, operation and maintenance costs. This is possible as we assume a necessity to sequester certain $CO_2$ amounts to reach climate neutrality.  The tractability also relies on using solutions, like the results from the perfect information model (M1), as starting points.

\section{Formal Model}
\label{sec:model}
In this section will describe the complete formal model. As shown in Figure \ref{fig:models} besides the main two stage regret model we have two auxiliary models. These are needed to compute best known solutions for certain scenario combinations and to evaluate the benefit of the regret model. The following index sets, parameters and variables are in common, but the auxiliary models omit some of them. 

\begin{tabular}{p{22mm}l}
	\textbf{Index Sets} & \\
	$\mathcal{V}$ & nodes \\
	$\mathcal{A}$ & arcs \\
	$\mathcal{S}$ & scenarios \\ 
	$\mathcal{C}$ & trends \\
	$\mathcal{O}$ & operations\\
\end{tabular}

\begin{tabular}{p{22mm}l}
	\textbf{{Parameters}} & \\
	$l_{ij}$ & length of an arc from node $i$ to $j$\\
	$m_{c}, y_{c}$ & cost for building a pipeline with trend $c$\\
	$Q_{c}^{max}$ & maximum diameter of pipeline with trend $c$\\
	$Q_{c}^{min}$ & minimum diameter of pipeline with trend $c$\\
	$d_{is}$ &  demand of node $i$ in scenario $s$ (negative for sinks)\\
	$d_{i1s}$ &  demand of node $i$ in scenario $s$ in period 1\\
	$n_1$ & years between first and second investment period\\
	$n_2$ & years between first investment period and end of planning horizon\\ 
	$OM$ &  factor for operation and maintenance\\
	$\tau$ &  discount rate\\
	$o_2^{max}$ & maximum pressure increase factor \\
	$o_2^{cost}$ & cost for pressure increase\\
\end{tabular}

\begin{tabular}{p{22mm}l}
	\textbf{Variables} &  \\
	$b_{ijc}$ & pipeline build on arc $(i,j)$ with trend $c$ in period 0\\
	$p_{ijc}$ & flow on arc $(i,j)$ with trend $c$ in period 0\\
	$b_{ijcs}$ & pipeline build on arc $(i,j)$ with trend $c$ in period 1 \\
	& considering scenario $s$\\
	$p_{ijcs}$ & flow on arc $(i,j)$ with trend $c$ in period 1 considering scenario $s$\\
	$f_{ijs}$ & total flow on arc $(i,j)$ in period 1 considering scenario $s$\\
	$r_{ijs}$ & restructure cost for arc $(i,j)$ in period 1 considering scenario $s$\\
	$u_{oijs}$ & use operation $o$ on arc $(i,j)$ in period $1$ considering scenario $s$\\
	$x$ & system regret
\end{tabular}

\setcounter{equation}{0}
\subsection{First stage}

In the first stage the model decides on the general layout of the network. It decides if to build a pipeline on an arc ($b_{ijc}$) and how much flow is put on that arc ($p_{ijc}$).

(\ref{eq:min_trend_s1}) to (\ref{eq:one_trend_s1}) ensure that the piece-wise linear trends are respected. This also includes that a flow is only allowed on an arc on which a pipeline is built. (\ref{eq:flow_s1}) ensures that all demand in the system is handled as well as flow is conserved throughout the pipeline network.
\begin{align}
	Q_{c}^{min} b_{ijc} &\leq p_{ijc} & (i,j) \in \mathcal{A}, c \in \mathcal{C} \label{eq:min_trend_s1}\\
	Q_{c}^{max} b_{ijc} &\geq  p_{ijc} &  (i,j) \in \mathcal{A}, c \in \mathcal{C} \label{eq:max_trend_s1}\\
	\sum_{c \in \mathcal{C}} b_{ijc} & \leq 1 &  (i,j) \in \mathcal{A}, c \in \mathcal{C} \label{eq:one_trend_s1}\\
	\sum_{(i,j) \in \mathcal{A}} \sum_{c \in \mathcal{C}} p_{ijc} - \sum_{(j,i) \in \mathcal{A}} \sum_{c \in \mathcal{C}} p_{jic}  &\geq d_{js}&  j \in \mathcal{V}, s=s_1 \label{eq:flow_s1}\\ 
	b_{ijc} &\in \{0,1\}& (i,j) \in \mathcal{A}, c \in \mathcal{C} \label{eq:build_binary_s1} \\
	p_{ijc} &\geq 0 &  (i,j) \in \mathcal{A}, c \in \mathcal{C} \label{eq:flow_no_neg_s1}
\end{align}

\subsection{Second stage}
In the second stage of the model all decision variables are indexed by $s$, so decisions are made for all possible realizations.

For each pipeline segment $(i,j)$ a binary decision $u_{oijs}$ is made if a certain upgrade function $o$ is used in scenario $s$. Constraint (\ref{eq:comb_options}) ensures that only one upgrade operation is used. Similar to (\ref{eq:flow_s1}), constraint (\ref{eq:comb_flow}) ensures that all demand is satisfied and flow conversation. A variable $f_{ijs}$ is introduced as this is possible to achieve by multiple different operations. This is manipulated by each operation with indicator constraints.

\begin{align}
	\sum\limits_{o \in O} u_{oijs}& = 1 &  (i,j) \in \mathcal{A},  s \in \mathcal{S} \label{eq:comb_options}\\
	\sum\limits_{(i,j) \in A} f_{ijs} - \sum\limits_{(j,i) \in \mathcal{A}} f_{jis} & \geq d_{j1s}&  j \in V,  s \in \mathcal{S} \label{eq:comb_flow}\\
	u_{oijs} &\in \{0,1\}& (i,j) \in \mathcal{A}, s \in \mathcal{S}, o \in O \label{eq:comb_options_binary}\\
	f_{ijs} &\geq 0 & (i,j) \in \mathcal{A}, s \in \mathcal{S} \label{eq:comb_flow_no_neg}
\end{align}

If option $o=1$ is selected (indicated by $u_{1ijs} = 1$), the model can build a parallel or a new line on an arc $(i,j)$. This is modeled analog to the first stage, but with scenarios and adding up both pipeline flows in (\ref{eq:indicator_loop_flow}). 

\begin{align}
	(u_{1ijs} == 1)  \Rightarrow f_{ijs} &= \sum\limits_{c \in \mathcal{C}} p_{ijc} + p_{ijcs} & (i,j) \in \mathcal{A},  s \in \mathcal{S} \label{eq:indicator_loop_flow}\\
	\sum_{c \in \mathcal{C}} b_{ijcs} & \leq 1 & (i,j) \in \mathcal{A},  s \in \mathcal{S}\\
	Q_{c}^{min} b_{ijcs} &\leq p_{ijcs} &  (i,j) \in \mathcal{A}, c \in \mathcal{C}, s \in \mathcal{S} \label{eq:min_trend_loop}\\
	Q_{c}^{max} b_{ijcs} &\geq  p_{ijcs} & (i,j) \in \mathcal{A}, c \in \mathcal{C}, s \in \mathcal{S} \label{eq:max_trend_loop}\\
	b_{ijcs} &\in \{0,1\}& (i,j) \in \mathcal{A}, c \in \mathcal{C}, s \in \mathcal{S} \label{eq:build_binary_loop} \\
	p_{ijcs} &\geq 0 &  (i,j) \in \mathcal{A}, c \in \mathcal{C}, s \in \mathcal{S} \label{eq:flow_no_neg_loop}
\end{align}
By selecting option $o=2$ (indicated by $u_{2ijs} = 1$) we allow to increase an already existing capacity by a factor $o_2^{max}$ by increasing the pressure in the pipeline (\ref{eq:indicator_pressure_flow}). Constraint (\ref{eq:indicator_pressure_restructure_cost}) enforces a restructuring cost $r_{ijs}$ on this.
\begin{align}
	(u_{2ijs} == 1)  \Rightarrow f_{ijs} &= \sum\limits_{c \in \mathcal{C}} p_{ijc} \cdot o_2^{max} &  (i,j) \in A,  s \in \mathcal{S} \label{eq:indicator_pressure_flow}\\
	(u_{2ijs} == 1)  \Rightarrow r_{ijs} &= I_{ij}^0 \cdot o_2^{cost} &  (i,j) \in A,  s \in \mathcal{S} \label{eq:indicator_pressure_restructure_cost}
\end{align}

To compute the system regret $x$, we need to calculate investment costs $I$, operating and maintenance costs $O$ and restructuring costs $R$ for each scenario. $I$ as well as $O$ occurring in $t=0$ are the same for all scenarios, but differ for $t=1$. As investments done in $t=1$ are not written-off completely at the end of the planning horizon, we need to calculate this rest value and subtract it from the $I^1$.  To compute the regret value, we further need the best solution for the case where we assume full knowledge which scenario occurs in $t=1$. To compute this best solution value $B_s$, we are using auxiliary model $M1$ (see section \ref{subsubsec:M1}). 

The target function (\ref{eq:target_regret}) in conjunction with constraint (\ref{eq:regret}) forces a \textit{min-max} regret optimization.
\begin{align}
	\text{min } & x \label{eq:target_regret}\\
	I^0 + O^0 + \frac{n_2 - n_1}{n_2} I_{s}^1 + O_{s}^1 + R_{s} - B_s &\leq x  &  s \in \mathcal{S}  \label{eq:regret}\\
	x & \geq 0& \label{eq:non_neg_regret}
\end{align}

with:
\begin{align*}
	I^0 &= \sum_{(i,j) \in A} \sum_{c \in \mathcal{C}} (p_{ijc0} \cdot m_{c} + b_{ijc0} \cdot y_c) \cdot l_{ij} \\
	O^0 &= \sum_{n=1}^{n_1} OM \cdot \dfrac{I^0}{(1+\tau)^n}\\
	I_{s}^1 &= \sum_{(i,j) \in A} \sum_{c \in \mathcal{C}} (p_{ijcs} \cdot m_{c} + b_{ijcs} \cdot y_c) \cdot l_{ij} \\
	R_{s} &= \sum_{(i,j) \in A} r_{ijs}\\
	O_{s}^1 &= \sum_{n=n_1}^{n_2} OM \cdot \dfrac{I^0 + I_{s}^1 + R_s}{(1+\tau)^n} \\
	B_s &= \text{Best solution for scenario $s$ with perfect information ($z_{M1}$)}
\end{align*}

\subsection{Auxiliary Models}
To calculate the best known solutions and also to compute initial solutions for the regret model, we derive the following models from the model described above.

\subsubsection{Perfect Information Model (M1)}
\label{subsubsec:M1}
The model M1 assumes full knowledge of scenario $s1$ and $s2$. It will make its investment decisions in $t = 0$ already considering all necessary upgrades for $t=1$. For the second stage, the set $S$ just contains one scenarios $s2$. The following objective function is used with $s = s2$, subject to (\ref{eq:min_trend_s1}) - (\ref{eq:indicator_pressure_restructure_cost}).
\begin{align}
	min\ z_{M1} &= I^0 + O^0 + \frac{n_2 - n_1}{n_2} I_{s}^1 + O_{s}^1 + R_{s}
\end{align}	

\subsubsection{Successive Information Model (M2)}
\label{subsubsec:M2}
The Model M2 works in two steps. It first computes the optimal solution network for one scenario $s1$ and for one period $t=0$. It assumes that the pipeline is used consistently for an amount $n_2$ of years. This problem resembles nearly the first stage decision problem in Section \ref{sec:model}. The following target function is minimizing investment cost $I^0$. This suffices, as operation and maintenance cost are directly depended on this and there are no further changes to the layout. The target function is subject to (\ref{eq:min_trend_s1}) - (\ref{eq:flow_no_neg_s1})
\begin{align}
	min\ z_{t0} &= \sum_{(i,j) \in \mathcal{A}} \sum _{c \in \mathcal{C}} (p_{ijt} \cdot m_{c} + b_{ijc} \cdot y_c) \cdot l_{ij}
\end{align}

In a second step, more information becomes available, namely the scenario $s2$ happening in $t=1$. The model now takes the resulting decision variables as parameters. Meaning all first stage variables ($b_{ijc}, p_{ijc} \ \ (i,j) \in A, c \in C$) become parameters. It than finds the optimal network for $t=1$ and the scenario $s2$. The following objective function is used, subject to (\ref{eq:comb_flow}) - (\ref{eq:indicator_pressure_restructure_cost}).
\begin{align}
	min\ z_{t1} &= \frac{n_2 - n_1}{n_2} I_{s}^1 + O_{s}^1 + R_{s}
\end{align}

\section{Case Study}
\label{sec:studies}

\begin{figure}[htbp]
	\centering
	\begin{subfigure}{.48\textwidth}
		\centering
		\includegraphics[width=0.9\linewidth]{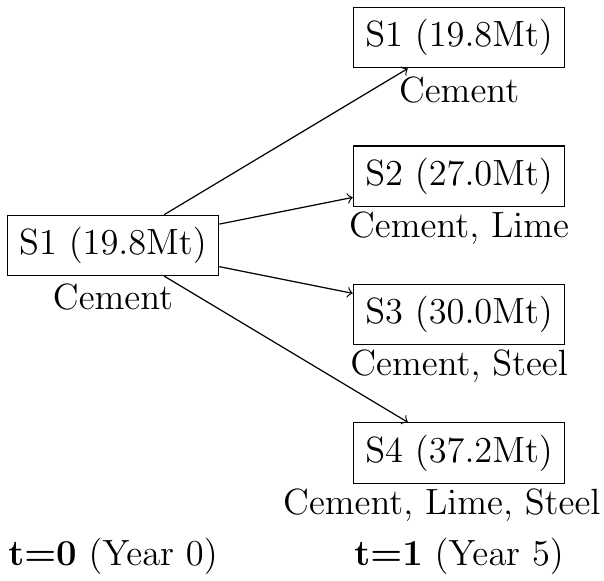}
		\caption{Scenario Tree}
		\label{subfig:scenario_tree}
	\end{subfigure}
	\begin{subfigure}{.48\textwidth}
		\centering
		\includegraphics[width=0.9\linewidth]{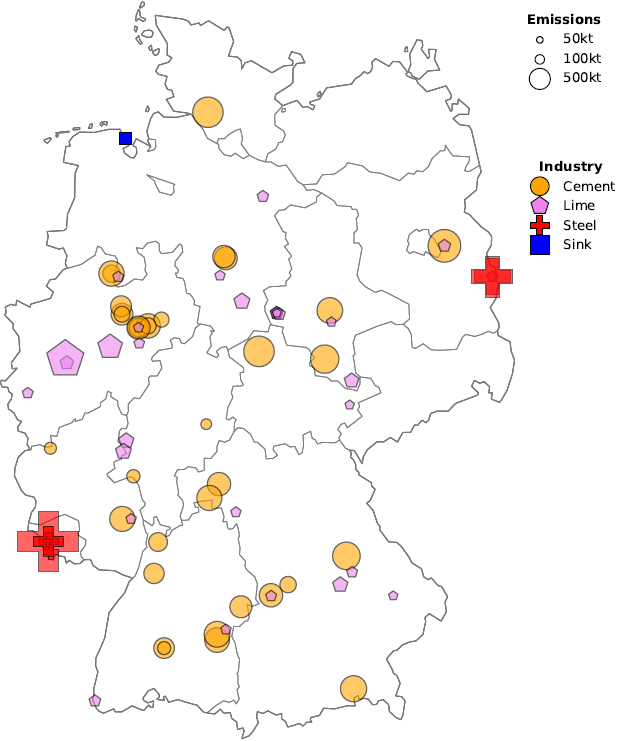}
		\caption{Geographical distribution of emitters}
		\label{subfig:geo_emitters}
	\end{subfigure}
	\caption{Scenarios}
	\label{fig:geo_emitters}
\end{figure}

To demonstrate our model, we present a case study on the development of a $CO_2$ pipeline network in Germany. The cement producers are assumed to be the first consumers of the initial network. This is attributed to the indispensability of CCS to decarbonize the sector due to the hard-to-abate process emissions \citep{iea_technology_2018}. Additionally, the sector's ongoing activities in planning of the relevant infrastructure demonstrates a clear and determined vision regarding the technology deployment \citep{vdz_anforderungen_2024}. Therefore,the uncertainty surrounding CCS demand is relatively minimal compared to other industrial sectors. Afterwards, four scenarios are considered for the network development (Figure \ref{subfig:scenario_tree}).

The first scenario (S1) assumes that no additional emitters are added to the network. The second (S2) and third (S3) scenarios assume that either the lime or steel emitters are added to the network, respectively. Finally, both the lime and steel emitters join the network in the fourth scenario (S4). An overview of the scenarios and their respective plant count and emissions can be found in Table \ref{tab:scenarios} \citep{dehst_greenhouse_2023}. The geographic distribution of the emitters is also depicted in (Figure \ref{subfig:geo_emitters}). The analyses consider the emitters with more than 100 kt $CO_2$ as well as the steel plants far away from the planned national hydrogen network. The study also considers a North Sea harbor (Wilhelmshaven) as a potential sink, based on \citet{vdz_anforderungen_2024}. All parameters can be found in \ref{app:case_params}. This also includes all the data sources for geographic features, population density, etc. as well as gas density, velocity, discount rate and cost parameters. The calculations are done with Gurobi version 11, on an Intel i7-13700K Processor with 16 cores with base speed of 3.4Ghz and 128GB of RAM.

\begin{table}[h]
	\centering
	\caption{Scenario details}
	\begin{tabular}{|c|c|c|c|}
		\hline
		Scenario & Description & Emitters & Emissions [Mt/a] \\
		\hline
		S1 & Cement  & 35 & 19.8 \\
		\hline
		S2 & Cement, Lime & 64 & 27.0 \\
		\hline
		S3 & Cement, Steel & 44 & 30.0 \\
		\hline
		S4 & Cement, Lime, Steel & 73 & 37.2 \\
		\hline
	\end{tabular}
	\label{tab:scenarios}
\end{table}

\begin{figure}[h]
	\centering
	\includegraphics[width=1.0\linewidth]{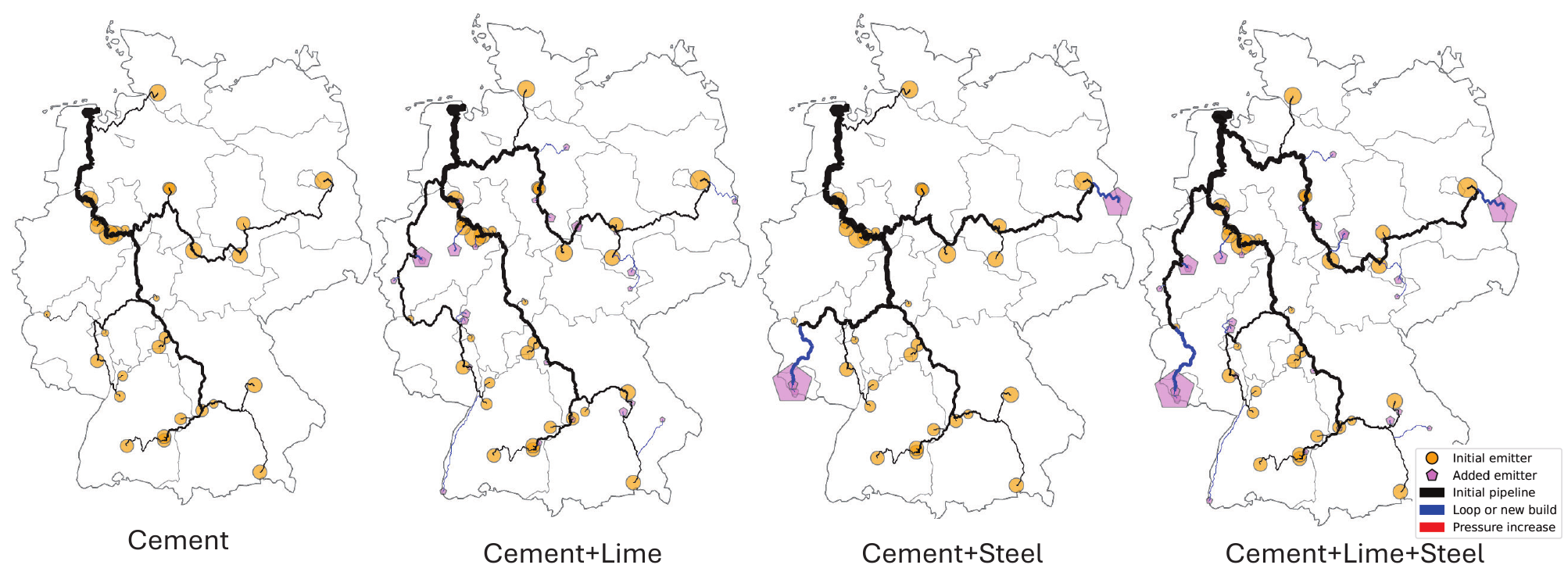}
	\caption{Perfect information model result layouts}
	\label{fig:case1fkt1}
\end{figure}

\begin{figure}[htbp]
	\centering
	\includegraphics[width=1.0\linewidth]{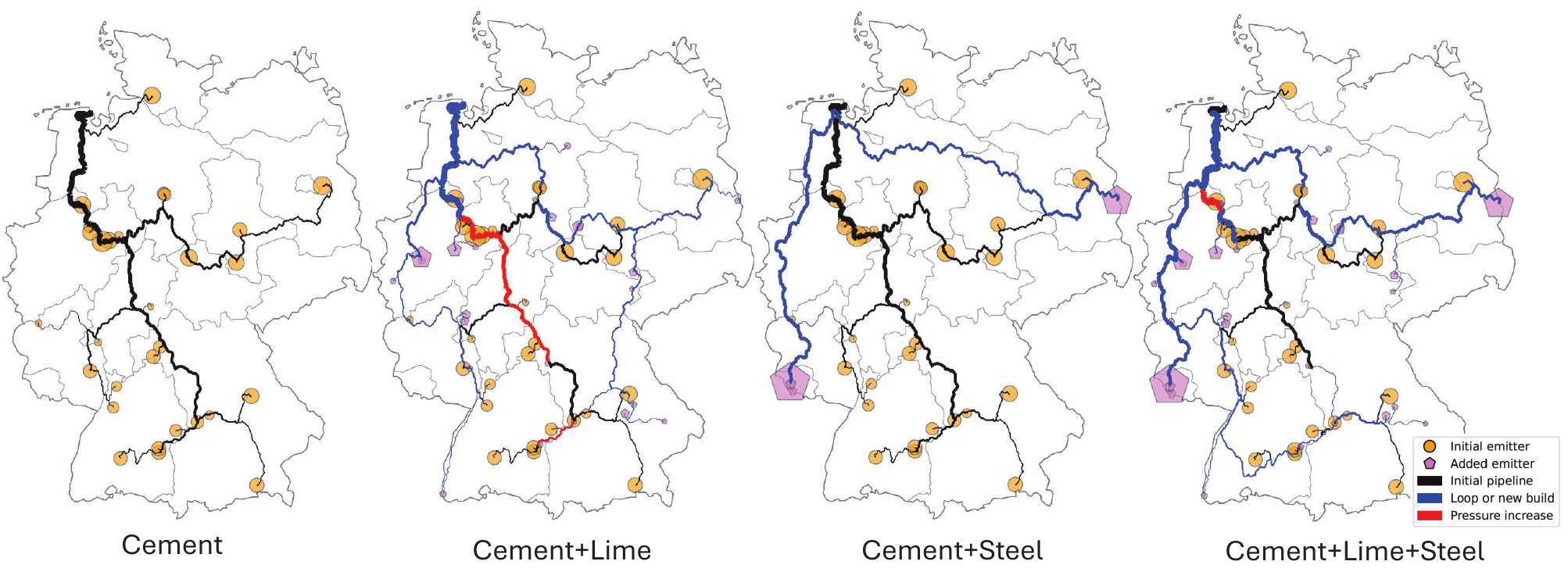}
	\caption{Successive information model result layouts}
	\label{fig:case1nkt1}
\end{figure}

\begin{figure}[h]
	\centering
	\includegraphics[width=1.0\linewidth]{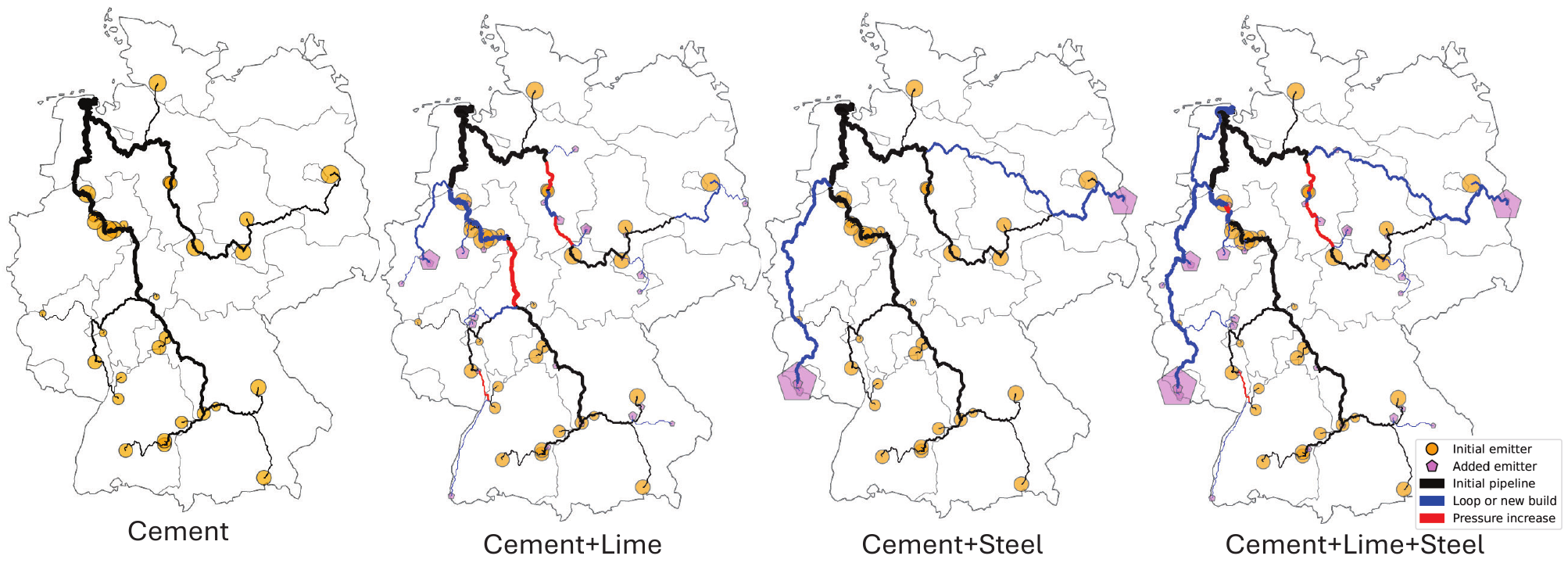}
	\caption{Regret model result layouts}
	\label{fig:case1regrett1}
\end{figure}

The results of the case study and the resulting layouts are provided in Table \ref{tab:system_cost_cs1} and \ref{tab:scenario_regret_and_benefit_cs1} as well as Figures \ref{fig:case1fkt1} - \ref{fig:case1regrett1}. Table \ref{tab:system_cost_cs1} presents the system costs for the three different models computed for the case study: the regret model, the perfect information model (M1), and the successive information model (M2). The perfect information model (M1) and the successive information model (M2) were run for each scenario combination separately. Table \ref{tab:scenario_regret_and_benefit_cs1} also shows the potential savings through perfect knowledge, the regret, and the benefit of the regret optimization. Herein, potential savings means the difference between the results of the perfect information model and the successive information model; benefit means the difference between the results of the regret model and the successive information model.

\begin{table}[h]
	\centering
	\caption{Total system costs}
	\begin{tabular}{|c|r|r|r|}
		\hline
		Scenarios & $z_{M1}$\tablefootnote{Perfect information model}[Mio €] & $z_{M2}$\tablefootnote{Successive information model}[Mio €] & $z_R$\tablefootnote{Regret model}[Mio €] \\
		\hline
		Cement (S1)& 3,141.325 & 3,141.325 & 3,497.485 \\
		\hline
		Cement, Lime (S2)& 3,979.308 &  4,997.963& 4,426.101 \\
		\hline
		Cement, Steel (S3)& 3,946.124 &  4,422.603& 4,400.765 \\
		\hline
		Cement, Lime,  Steel (S4)& 4,702.754 & 5,653.400 & 5,130.106 \\
		\hline
	\end{tabular}
	\label{tab:system_cost_cs1}
\end{table} 

\begin{table}[h]
	\centering
	\caption{Potential savings, regret and benefit}
	\begin{tabular}{|c|r|r|r|}
		\hline
		Scenarios & Potential[Mio €]
		& Regret[Mio €] & Benefit[Mio €] \\
		\hline
		Cement  (S1)& 0.000 & 356.161 & -356.161 \\
		\hline
		Cement, Lime (S2)& 1,018.655 & 446.792 & 571.863 \\
		\hline
		Cement, Steel (S3)& 476.479 &  454.640 & 21.839 \\
		\hline
		Cement, Lime, Steel (S4)& 997.411 & 427.352 & 523.295 \\
		\hline
	\end{tabular}
	
	\label{tab:scenario_regret_and_benefit_cs1}
\end{table}

The layouts in Figures \ref{fig:case1fkt1} - \ref{fig:case1regrett1} are visualized to show both investment periods in the same map. For t = 0, circles represent the emitters and black lines show the pipeline build. For the second period, the pentagons indicate the added emitters. The other lines either illustrate the construction o  a new pipeline or a parallel line (upgrade option 1) or show that the pressure in the pipeline is increased (upgrade option 2).

For the perfect information model (M1), the costs of the different scenarios range between 3,141 Mio. and 4,702 Mio. €. The outcomes of the perfect information model are characterized by oversizing to accommodate all requested demand. Subsequently, only the connecting pipelines to the existing network are built in t = 1. As shown in Figure \ref{fig:case1fkt1}, no parallel lines or pressure increases were used. Using the successive information model (M2), the initial network of scenario S1 is constructed first and then upgraded at t =1. Therefore, as expected, this approach exhibits significant increases in cost. In the cement and lime scenario (S2), the cost increase from 3,979 Mio. € to 4,998 Mio. € (+1,019 Mio. €). For the cement, lime and steel scenario (S4), there is an analogous increase of 997 Mio. € from 4,703 Mio. € to 5,653 Mio. €. In the cement and steel scenario (S3), there is also an increase of 476 Mio. € , which is significantly lower compared to the other scenarios.

Figure \ref{fig:case1nkt1} depicts the layout changes of the successive information model. Compared to Figure \ref{fig:case1fkt1}, significant interventions all over the network are observed. Herein, a mixture out of pressure increases, parallel lines and new pipelines are used in S2. Besides increasing the pressure in the main trunk of the pipeline, two smaller networks in the east and west are installed. Contrariwise, the initial layout is ignored in S3, with new lines being built primarily in areas where no existing pipeline is present. For S4, nearly no pressure increase is used. One new eastern trunk is established, which also connects to sources in Bavaria. From east to west, the new trunk is mostly following the initial layouts, but sometimes establishing new lines in not used before areas.

Using the regret model, the scenarios S2, S3 and S4 perform better than with the successive information model, with S2 and S4 having improvements of over 500 Mio. €. S3 has a small benefit of around 22 Mio. € and S1 will lead to extra cost of 356 Mio. €. However, the biggest regret has the scenario S3 with 455 Mio. €. The network presented in Figure \ref{fig:case1regrett1} showcases a rather clean layout. In many sections, the network is just oversized to allow simple connections in t = 1. However, the network makes use of pressure increases in specific parts of the network (S2, S3). The network in t = 0 looks very similar to the perfect knowledge model for S4. With the notable difference, that a main trunk for the steel plant to Saarland (south-west) is not built.

\section{Conclusion}
\label{sec:conclusion}
We propose a novel regret minimizing optimization model for CCS pipeline network. In a case study we prove that our model is able to identify networks that are adaptable to changing scenarios over time. It computes the trade-offs of these adapted and oversize layouts. Thus, our framework allows decision makers to access and quantify these these trade-offs and finally make good decisions.

The case study results demonstrate that oversizing the pipeline is clearly the preferred option if perfect information is available. Contrariwise, building without considering potential future scenarios is not a viable option, as it proves to be highly inefficient. Indeed, both models represent two extreme cases. On the one hand, it is strategically naive to overlook the potential future changes. That is why the initial layout of the successive information model has witnessed significant changes to allow the integration of additional emitters. On the other hand, the potential changes are also associated with several uncertainties, which is why the perfect knowledge model is not effective in the real applications. Herein, the regret optimization model emerges as a viable approach to tackle this dilemma. As the results of regret optimization model demonstrate, the required alterations are significantly lower than with the successive information model (M2). It should be highlighted that the regret model will still incur alterations as it is impossible to eliminate uncertainties. Nonetheless, while there may not be a solution without regret, the model can serve as the optimal approach to minimize it.

Looking ahead, the model can undergo upgrades in different directions. Firstly, exploring multi-stage or extended investment periods could provide a better understanding of long-term dynamics and potential outcomes. Secondly, considering reduction scenarios into our framework can also be a useful integration. While the model can theoretically handle decreasing values, sensible aggregated scenario rules would need to be determined as well as downgrade options for the pipelines. Thirdly, investigating alternative transportation connections, such as rail, could offer insights into diversifying infrastructure options and mitigating risks associated with single-mode dependencies. Lastly, increasing the model’s scalability can enhance the effectiveness of our approach. While performance is satisfactory at present, the model's capacity reaches its limit with industry-scale cases. Hence, improving scalability to incorporate additional cases, scenarios and data inputs will be crucial for sustaining the model's utility and relevance in increasingly complex decision-making environments.

\section*{Declaration of competing interest}
The authors report there are no competing interests to declare.

\section*{CRediT }
\textbf{Stephan Bogs:} Conceptualization, Methodology, Software, Formal analysis, Investigation, Writing - Original draft, Visualization. \textbf{Ali Abdelshafy:} Conceptualization, Validation, Investigation, Data curation, Writing - Original draft, Writing - Review and editing, Visualization. \textbf{Grit Walther:} Conceptualization, Resources, Writing - Review and editing, Supervision, Project administration, Funding acquisition.

\appendix

\section{Graph from geographic data}
\label{app:graph}
With the following steps we derive the candidate graph from the geographic data and emitters
\begin{enumerate}
	\item Split the geographic area in $1.5km \cdot 1,5km$ squares. 
	\item Apply multipliers for: Population density, slope, CDDA Zones and nature parks, lakes and waterways, infrastructure like existing pipelines, motorways and railroads.
	\item Inject emitter and sinks into the grid, this includes also emitter that are not present in the all scenarios to ensure a consistent graphs between the time periods. All emitter that are in the same square are treated as one node, with their emissions summed up.
	\item Use Delaunay triangulation between the all emitters (including hypothetical ones) and sinks.
	\item Calculate the shortest path between the triangulated points. For this we consider each square a node in a graph. Each orthogonal neighbor is connected with an edge that has the weight according to the multiplier. Each diagonal neighbor has an arc that has the weight multiplied with $\sqrt{2}$.
	\item We save the shortest paths with theirs nodes into a new undirected graph.
	\item In the new graph merge all nodes with degree of two thus only remaining nodes with degree greater equal three (transport nodes or emitters, sinks) or one (just emitters and sinks).
	\item Convert the undirected graph into a directed graph.		
\end{enumerate}

\section{Case Study Parameters}
\label{app:case_params}
For nearly all parameters we follow \citet{yeates_industrial_2024}. Gas Density of $p = 900\frac{kg}{m^3}$ and a velocity of $v = 3\frac{m}{s}$. We use \cite{parker_using_2003} formula adjusted by inflation to 2022 and a discount rate of $5\%$. First decision are made have $n_1=5$ years. We assume a $n_2=25$ year operational period which is the official write-off period for pipelines in Germany \citep{bundesfinanzministerium_afa-tabelle_1995}

For the generation of the graph we use the data listed in Table \ref{tab:datasets} and the multipliers listed in Table \ref{tab:multipliers}.

\begin{table}[htbp]
	\centering
	\caption{Multipliers used}	
	\label{tab:multipliers}		
	\begin{tabular}{|l|c|}
		\hline
		\textbf{Criterion} & \textbf{Multiplier} \\
		\hline
		Population density/$km^2$ ($<$250) & 1 \\
		\hline
		Population density/$km^2$  (250-500) & 4 \\
		\hline
		Population density/$km^2$ (500-2000) & 9 \\
		\hline
		Population density/$km^2$  (2000-4000) & 16 \\
		\hline
		Population density/$km^2$  (4000-8000) & 25 \\
		\hline
		Population density/$km^2$  ($>$8000) & 36 \\
		\hline
		Pre-existing pipelines & 0.25 \\
		\hline
		Railroads & 3 \\
		\hline
		Motorways & 3 \\
		\hline
		Rivers, lakes, and transitional waters & 10 \\
		\hline
		CDDA protected areas (excl. National parks) & 10 \\
		\hline
		National parks & 30 \\
		\hline
		Terrain slope [0°-90°] & 1-20 \\
		\hline
	\end{tabular}
\end{table}

\begin{table}[htbp]
	\centering
	\caption{Datasets used}	
	\label{tab:datasets}			
	\begin{tabular}{|l|l|}
		\hline
		\textbf{Data} & \textbf{Source} \\
		\hline
		Geographic outlines &  \citet{eurostat_nuts_2021} \\
		\hline
		Pre-existing pipelines &  \citet{kunz_electricity_2017} \\
		\hline
		Population density&  \citet{batista_e_silva_f_jrc-geostat_2021}\\
		\hline
		CDDA protected areas and national parks &  \citet{european_environment_agency_cdda_2020} \\
		\hline
		Road and rail infrastructure & \citet{kelso_natural-earth-vector_2024} \\
		\hline
		Topography &  \citet{copernicus_copernicus_2022} \\
		\hline
		Waterways, lakes, transitional waters &  \citet{kelso_natural-earth-vector_2024} \\
		\hline
	\end{tabular}
\end{table}

\bibliographystyle{elsarticle-harv}
\bibliography{sources}
	
\end{document}